\documentclass{amsart}
\usepackage{url}
\usepackage{xcolor}

\newtheorem{theorem}{Theorem}
\newtheorem{lemma}{Lemma}

\newcommand{\norm}[1]{\lVert#1\rVert}

\newcommand{\beq}[1]{\begin{equation}\label{#1}}
\newcommand{\eeq}{\end{equation}}

\newcommand{\beqan}{\begin{eqnarray*}}
\newcommand{\eeqan}{\end{eqnarray*}}

\numberwithin{equation}{section}
\numberwithin{lemma}{section}

\title{The Perron-Frobenius Theorem for Markov Semigroups}

\author{Omar Hijab}
\address{Department of Mathematics, Temple University, Philadelphia, PA 19122}
\email{hijab@temple.edu}
 
 \keywords{Markov semigroup, Perron-Frobenius, Donsker-Varadhan}
 \subjclass[2000]{47D07,58J65,52A40}
\date{\today}

\begin{document}

\begin{abstract}
Let $P^V_t$, $t\ge0$, be the Schrodinger semigroup associated to a potential
$V$ and Markov semigroup $P_t$, $t\ge0$, on $C(X)$. Existence is established of
a left eigenvector and right eigenvector corresponding to the spectral radius
$e^{\lambda_0t}$ of $P^V_t$, simultaneously for all $t\ge0$.
This is derived with no compactness assumption on the semigroup operators.
\end{abstract}

\maketitle

%\tableofcontents

\section{Introduction}

Let $X$ be a compact metric space and let $P_t$, $t\ge0$, be a Markov semigroup on $C(X)$ with 
generator $L$.  Given $V$ in $C(X)$ let $P_t^V$, $t\ge0$, denote the Schrodinger semigroup on $C(X)$ generated by $L+V$. 
Then the principal eigenvalue
$$\lambda_0(V)\equiv\lim_{t\to\infty} \frac1t\log\norm{P^V_t}$$
is given by the Donsker-Varadhan formula \cite{DV}
$$\lambda_0(V)= \sup_\mu \left(\int_XV\,d\mu-I(\mu)\right),$$
where the supremum is over probability measures $\mu$ on $X$, and 
$$I(\mu)\equiv-\inf_{u\in\mathcal D^+}\int_X\frac{Lu}{u}\,d\mu.$$
Here the infimum is over positive $u$ in the domain $\mathcal D$ of $L$.

An {\em equilibrium measure} is a measure $\mu$ achieving the supremum in the Donsker-Varadhan formula. 
Let $\lambda_0=\lambda_0(V)$.

A {\em ground state relative to $\mu$} is a Borel function $\psi$ satisfying $\psi>0$ a.s. $\mu$ and
$$e^{-\lambda_0t}P^V_t\psi=\psi,\qquad a.s. \mu, t\ge0.$$

A {\em ground measure} is a measure $\pi$ satisfying
\beq{eq:Pinv}
\int_Xe^{-\lambda_0t}P^V_tf\,d\pi=\int_Xf\,d\pi,\qquad t\ge0
\eeq
for $f$ in $C(X)$.

\begin{theorem} 
\label{theorem:triple}
Suppose $\pi$ and $\mu$ are measures with $\mu<<\pi$.
Suppose also $\psi=d\mu/d\pi$ satisfies $\psi\log\psi\in L^1(\pi)$.  Then the following hold.
\begin{itemize}
\item If $\pi$ is a ground measure and $\psi$ is a ground state relative to $\mu$,  then $\mu$ is an equilibrium measure.
\item If $\pi$ is a ground measure and $\mu$ is an equilibrium measure,  then $\psi$ is a ground state relative to $\mu$.
\item If $\mu$ is an equilibrium measure and $\psi$ is a ground state relative to $\mu$,  then $\pi$ is a ground measure.
\end{itemize}
\end{theorem}

Here is the Perron-Frobenius Theorem in this setting.

\begin{theorem}
\label{theorem:pf}
Fix $V$ in $C(X)$, suppose 
$$e^{-\lambda_0t}\norm{P^V_t}\le C,\qquad t\ge0,$$
for some $C>0$, and let $\mu$ be an equilibrium measure. Then there is a ground measure $\pi$ satisfying $\mu<<\pi$,
$\psi=d\mu/d\pi$ is a ground state relative to $\mu$, and $\psi\log\psi\in L^1(\pi)$.
\end{theorem}

This says there is a nonnegative right eigenvector $\psi$ and a nonnegative left eigenvector $\pi$, corresponding to the spectral radius $e^{\lambda_0t}$ of $P^V_t$, simultaneously for all $t\ge0$. 

Neither the hypothesis nor the conclusion hold when $L+V$ is a Jordan block: Take $X=\{0,1\}$ and 
$$L=\begin{pmatrix} 0 & 0 \\ 1 &  -1 \end{pmatrix},\qquad V=\begin{pmatrix} 0 \\ 1 \end{pmatrix}.$$
In this case every measure $\mu=(p,1-p)$ is an equilibrium measure and there is a unique ground measure $\pi=(0,1)$. 

In this generality, there is no guarantee of uniqueness of $\mu$, $\pi$ or $\psi$. 

If $P_t$, $t\ge0$, is self-adjoint in $L^2(\rho)$ for some measure $\rho$ on $X$, then, under suitable conditions, the Donsker-Varadhan formula reduces \cite{DV1} to the classical Rayleigh-Ritz formula for the principal eigenvalue, and  every ground state $\psi$ relative to $\rho$ yields a ground measure $d\pi=\psi\,d\rho$ and an equilibrium measure $d\mu=\psi\,d\pi=\psi^2\,d\rho$. In the self-adjoint case, existence of ground states is classical (Gross \cite{G}).

The Perron-Frobenius Theorem (with uniqueness) is known to hold in $L^p$ for positive operators in these cases: for finite irreducible matrices (Friedland \cite{F}, Friedland-Karlin \cite{FK}, Sternberg \cite{S}), under a positivity improving property (Aida \cite{A}), under a spectral gap condition (Gong-Wu \cite{GW}), and under uniform integrability or irreducibility (Wu \cite{W}). 

The existence of ground measures for positive operators and the existence of ground states for compact positive operators are classical results due to Krein-Rutman \cite{KR}, see also Schaefer \cite{HHS}, \cite{HHS1}. In general, ground states do not exist pointwise everywhere on $X$, for example when $L\equiv0$. The novelty of the above result is the handling of the general non-compact case by interpreting ground states as densities against ground measures, and the link with equilibrium measures.

The techniques used here are self-contained and follow  the original papers \cite{DV}, \cite{DV1}.   Preliminaries are discussed in section \ref{sec:pre}, equilibrium measures in section \ref{sec:em}, ground measures and ground states in section \ref{sec:gsgm}, and entropy in section \ref{sec:ent}.  The theorems are proved in section \ref{sec:thm}.

\section{Preliminaries}
\label{sec:pre}

Let $X$ be a compact metric space, let $C(X)$ denote the space of real continuous functions with the sup
norm $\norm{\cdot}$, and let $M(X)$ denote the space of Borel probability measures with the topology of weak convergence. Then $M(X)$
is a compact metric space. Throughout $\mu(f)$ denotes the integral of $f$ against $\mu$ and all measures are probability measures.

A {\em Markov semigroup on $C(X)$} is a strongly continuous semigroup $P_t:C(X)\to C(X)$, $t\ge0$, 
preserving positivity, $P_tf\ge0$, for $f\ge0$, and satisfying $P_t1=1$. 

The subspace $\mathcal D\subset C(X)$ of functions $f\in C(X)$
for which the limit
\begin{equation}
\label{eq:gen}
\left.\frac{d}{dt}\right|_{t=0} P_tf
\end{equation}
exists in $C(X)$ is dense. If $Lf$ is defined to be this limit, then the operator $L$ is
the {\em generator} of $P_t$, $t\ge0$, on $C(X)$.

Given $V$ in $C(X)$, the {\em Schrodinger semigroup on $C(X)$} associated to $V$ is the 
unique strongly continuous semigroup $P^V_t:C(X)\to C(X)$, $t\ge0$, preserving positivity, $P^V_tf\ge0$, for $f\ge0$,
with generator $L+V$, in the sense the limit
\beq{eq:pv}
\left.\frac{d}{dt}\right|_{t=0} P^V_tf
\eeq
exists in $C(X)$ iff $f\in\mathcal D$, and equals $Lf+Vf$. The Schrodinger semigroup may be constructed
as the unique solution of 
\beq{eq:schro}
P^V_tf=P_tf+\int_0^t P_{t-s}VP^V_sf\,ds,\qquad t\ge0.
\eeq
for $f\in C(X)$.  For $f\ge0$, (\ref{eq:schro}) implies
\beq{eq:major}
e^{\min Vt} P_tf\le P^V_tf\le e^{\max Vt}P_tf,\qquad t\ge0.
\eeq
Since $\norm{P^V_t}=\max P^V_t1$, this implies 
$$\min V\le \lambda_0(V)\le \max V.$$

For $\mu$ in $M(X)$ and $\lambda_0=\lambda_0(V)$, let
$$I^V(\mu)\equiv I(\mu)-\int_XV\,d\mu+\lambda_0=-\inf_{u\in\mathcal D^+}\int_X\frac{(L+V-\lambda_0)u}{u}\,d\mu$$
Then $I^0(\mu)=I(\mu)$ and  $I^V(\mu)=0$ iff $\mu$ is an equilibrium measure for $V$.

\begin{lemma}
For $V$ in $C(X)$, $I^V$ is lower semicontinuous, convex, and nonnegative. 
In particular, $I$ is lower semicontinuous, convex, and nonnegative.
\end{lemma}

\proof
Lower semicontinuity and convexity follow from the fact that $I^V$ is the supremum of continuous affine functions.
The Donsker-Varadhan formula implies $I^V$ is nonnegative.\qed

Let $L^1(\mu)$ denote the $\mu$-integrable Borel functions on $X$ with norm 
$$\norm{f}_{L^1(\mu)}=\int_X|f|\,d\mu=\mu(|f|).$$

For $V$ in $C(X)$, and let $P^V_t$, $t\ge0$, be the Schrodinger semigroup on $C(X)$. Let $\mathcal D$ be the domain of the generator $L$ of $P_t$, $t\ge0$, $\mathcal D^+$ the positive functions in $\mathcal D$, and $C^+(X)$ the positive functions in $C(X)$. 

By positivity, there is a family $(t,x)\mapsto p^V_t(x,\cdot)$ of Borel measures on $X$ such that the Schrodinger semigroup may be written
\beq{eq:posmeas}
P_t^Vf(x)=\int_Xf(y)\,p^V_t(x,dy)
\eeq
for $t\ge0$, $x\in X$, and $f\in C(X)$. Hence $P^V_tf(x)\le+\infty$  is defined for nonnegative Borel $f$ for all $t\ge0$ and $x\in X$.

Let $B(X)$ denote the bounded Borel functions on $X$ and let $B^+(X)$ denote the positive Borel functions $u$ on $X$ with $f=\log u\in B(X)$. We say $f_n\in B(X)$  converges {\em boundedly} to $f\in B(X)$ if $f_n\to f$ pointwise everywhere and there is a $C>0$ with $|f_n|\le C$, $n\ge1$. 

\begin{lemma}
\label{lemma:boundedly}
If $f_n\to f$ boundedly then $u_n=e^{f_n}\to u=e^f$ boundedly, $P^V_tu_n\to P^V_tu$ boundedly, and 
$$
\log\left(\frac{e^{-\lambda_0t}P^{V}_tu_n}{u_n}\right)\to \log\left(\frac{e^{-\lambda_0t}P^{V}_tu}{u}\right)
$$
boundedly.
\end{lemma}

\proof 
If $|f_n|\le C$, $n\ge1$, then $u_n\to u$ pointwise. Since $e^{-C}\le u_n\le e^C$, $u_n\to u$ boundedly. Hence $P^V_tu_n\to P^V_tu$ pointwise by the dominated convergence theorem. By  (\ref{eq:major}), $P^V_tu_n\to P^V_tu$ boundedly. Since
$$
\left|\log\left(\frac{e^{-\lambda_0t}P^{V}_tu}{u}\right)\right|\le t(\max V-\min V)+\log\left(\frac{\sup u}{\inf u}\right),
$$
the last statement follows.\qed

\section{Equilibrium Measures}
\label{sec:em}

\begin{lemma}For $u$ in $B^+(X)$,
\label{lemma:Iapprox}
\beq{eq:log}
\int_X\log\left(\frac{e^{-\lambda_0t}P^{V}_tu}{u}\right)\,d\mu\ge-tI^V(\mu),\qquad t\ge0.
\eeq
\end{lemma}

\proof
By definition of $I^V(\mu)$,
\beq{eq:equil}
\int_X\frac{(L+V-\lambda_0)u}{u}\,d\mu\ge-I^V(\mu),\qquad u\in\mathcal D^+.
\eeq
For $t=0$, (\ref{eq:log}) is an equality. Moreover for $t>0$ and $u\in\mathcal D^+$, we
have $e^{-\lambda_0t}P^V_tu\in\mathcal D^+$ and
$$\frac d{dt} \int_X\log\left(\frac{e^{-\lambda_0t}P^{V}_tu}{u}\right)\,d\mu=\int_X\frac{(L+V-\lambda_0)(e^{-\lambda_0t}P^V_tu)}{e^{-\lambda_0t}P^V_tu}\,d\mu\ge -I^V(\mu).$$
This establishes (\ref{eq:log}) for $u\in\mathcal D^+$. Since $\mathcal D^+$ is dense in $C^+(X)$, (\ref{eq:log}) is valid for $u$ in $C^+(X)$. 

Now if $\log u_n\to \log u$ boundedly and (\ref{eq:log}) holds for $u_n$, $n\ge1$,  then by Lemma \ref{lemma:boundedly}, (\ref{eq:log}) holds for $u$.
Thus the class of Borel functions $f=\log u$ for which (\ref{eq:log})  holds is closed under bounded convergence. Thus (\ref{eq:log}) holds for all $u\in B^+(X)$.
\qed

The following strengthening of Lemma \ref{lemma:Iapprox}  is necessary below. 

\begin{lemma} 
\label{lemma:strong}
Let $u>0$ Borel satisfy $\log u\in L^1(\mu)$. Then for $t\ge0$,
\beq{eq:logP}
tI^V(\mu)+\int_X\log^+\left(\frac{e^{-\lambda_0t}P^{V}_tu}{u}\right)\,d\mu\ge\int_X\log^-\left(\frac{e^{-\lambda_0t}P^{V}_tu}{u}\right)\,d\mu.
\eeq
Here the integrals may be infinite.
\end{lemma}

\proof Without loss of generality, assume $I^V(\mu)<\infty$. 

Assume in addition $u\ge\delta>0$ and let $u_n=u\wedge n$, $n\ge1$. Then $u_n$ is in $B^+(X)$, (\ref{eq:log}) holds with $u_n$, and $u\ge u_n$,  hence
$$\int_X \log\left(\frac{e^{-\lambda_0t}P^V_tu}{u}\right)\,d\mu\ge \int_X \log\left(\frac{e^{-\lambda_0t}P^V_tu_n}{u}\right)\,d\mu
\ge - \int_{u>n}\log u\,d\mu-tI^V(\mu).$$
Letting $n\to\infty$ yields  (\ref{eq:log}) hence (\ref{eq:logP}) for $\log u$ in $L^1(\mu)$, provided $u\ge\delta>0$. Here the left side of (\ref{eq:logP}) may be infinite, but the right side is finite.

Now for $u>0$ Borel with $\log u$ in $L^1(\mu)$, let  $u_\delta=u\vee\delta$. Then by what we just derived,
$$
tI^V(\mu)+\int_X\log^+\left(\frac{e^{-\lambda_0t}P^V_tu_\delta }{u_\delta}\right)\,d\mu\ge\int_X\log^-\left(\frac{e^{-\lambda_0t}P^V_tu_\delta}{u_\delta}\right)\,d\mu
$$
so
\beq{eq:logPdelta}
tI^V(\mu).+\int_X\log^+\left(\frac{e^{-\lambda_0t}P^V_tu_\delta }{u}\right)\,d\mu
\ge\int_X\log^-\left(\frac{e^{-\lambda_0t}P^V_tu_\delta}{u}\right)\,d\mu+\int_{u<\delta}\log u\,d\mu.
\eeq
We may assume 
$$\int_X\log^+\left(\frac{e^{-\lambda_0t}P^V_tu}{u}\right)\,d\mu<\infty,$$
otherwise (\ref{eq:logP}) is vacuously true. 
To establish (\ref{eq:logP}), we pass to the limit $\delta\downarrow0$ in (\ref{eq:logPdelta}). Since $\log u\in L^1(\mu)$ and
$$\log^-\left(\frac{e^{-\lambda_0t}P^V_tu_\delta}{u}\right),\qquad \delta>0,$$
is an increasing sequence in $\delta>0$, the right side of (\ref{eq:logPdelta}) converges to the right side of (\ref{eq:logP}) as $\delta\downarrow0$.
Since $u_\delta\le u+\delta$, 
$$\log^+\left(\frac{e^{-\lambda_0t}P^V_tu_\delta }{u}\right)\le \log^+\left(\frac{e^{-\lambda_0t}P^V_tu }{u}\right)+\log^-u+C,\qquad \delta\le1,$$
hence the dominated convergence theorem applies. Thus the left side of (\ref{eq:logPdelta}) converges to the left side of (\ref{eq:logP}) as
$\delta\downarrow0$.
\qed

If $\psi>0$ is Borel, for $u\ge0$ Borel, define
\beq{eq:Ppsi}
P^{V,\psi}_tu\equiv e^{-\lambda_0t}\frac{P^V_t(\psi u)}{\psi},\qquad t\ge0.
\eeq

\begin{lemma}
\label{lemma:ineqQ}
Suppose $\log\psi\in L^1(\mu)$ and let $u>0$ Borel satisfy $\log u\in L^1(\mu)$. Then for $t\ge0$,
\beq{eq:logQ}
tI^V(\mu)+\int_X\log^+\left(\frac{P^{V,\psi}_tu}{u}\right)\,d\mu\ge\int_X\log^-\left(\frac{P^{V,\psi}_tu}{u}\right)\,d\mu.
\eeq
Here the integrals may be infinite.
\end{lemma}

\proof Since $\log\psi$ is in  $L^1(\mu)$, $\log(u\psi)$ is in $L^1(\mu)$ iff $\log u$ is in $L^1(\mu)$. Now apply Lemma
\ref{lemma:strong}.\qed

\begin{lemma}
\label{lemma:nonneg}
Suppose $\log\psi\in L^1(\mu)$. Then $\mu$ is an equilibrium measure iff for all $u>0$ Borel satisfying $\log u\in L^1(\mu)$,
 $$
\int_X\log^+\left(\frac{P^{V,\psi}_tu}{u}\right)\,d\mu\ge\int_X\log^-\left(\frac{P^{V,\psi}_tu}{u}\right)\,d\mu,\qquad t\ge0.
$$
\end{lemma}

\proof
If $\mu$ is an equilibrium measure, $I^V(\mu)=0$ so the result follows from Lemma \ref{lemma:ineqQ}.
Conversely, assume the inequality holds for all $u>0$ satisfying $\log u\in L^1(\mu)$.
For $u\in C^+(X)$, the function $u/\psi$  satisfies $\log (u/\psi)\in L^1(\mu)$. Inserting $u/\psi$ in the inequality yields
 $$
 \int_X\log^+\left(\frac{e^{-\lambda_0t}P^V_tu}{u}\right)\,d\mu\ge\int_X\log^-\left(\frac{e^{-\lambda_0t}P^V_tu}{u}\right)\,d\mu.
$$
For $u$ in $C^+(X)$, these integrals are finite hence
$$
 \int_X\log\left(\frac{e^{-\lambda_0t}P^V_tu}{u}\right)\,d\mu\ge0.
 $$
Now for $u\in\mathcal D^+$,
$$e^{-\lambda_0t}P^V_tu=u+t(L+V-\lambda_0)u+o(t),\qquad t\to0,$$
uniformly on $X$. Since $u>0$,
$$\frac{e^{-\lambda_0t}P^V_tu}{u}=1+t\frac{(L+V-\lambda_0)u}{u}+o(t),\qquad t\to0,$$
uniformly on $X$. Thus
$$\log\left(\frac{e^{-\lambda_0t}P^V_tu}{u}\right)=t\frac{(L+V-\lambda_0)u}{u}+o(t),\qquad t\to0,$$
uniformly on $X$, hence
$$
\lim_{t\to0}\frac1t \int_X \log\left(\frac{e^{-\lambda_0t}P^V_tu}{u}\right)\,d\mu=\int_X\frac{(L+V-\lambda_0)u}{u}\,d\mu.
$$
This implies
$$\int_X\frac{(L+V-\lambda_0)u}{u}\,d\mu\ge0. 
$$
This implies $I^V(\mu)\le0$, hence $I^V(\mu)=0$.\qed

\section{Ground Measures and Ground States}
\label{sec:gsgm}

\begin{lemma} 
\label{lemma:PL}
If $\pi$ is a ground measure, then
$e^{-\lambda_0t}P^V_t$, $t\ge0$, is a strongly continuous contraction semigroup on $L^1(\pi)$ and (\ref{eq:Pinv}) holds for $f$ in $L^1(\pi)$.
\end{lemma}

\proof  (\ref{eq:Pinv}) implies
$$\norm{e^{-\lambda_0t}P^V_tf}_{L^1(\pi)}=\int_X|e^{-\lambda_0t}P^V_tf|\,d\pi\le \int_Xe^{-\lambda_0t}P^V_t|f|\,d\pi =\int_X|f|\,d\pi=\norm{f}_{L^1(\pi)}$$
for $f$ in $C(X)$. Let $f$ be in $L^1(\pi)$ and choose $f_n\in C(X)$ converging to $f$ in $L^1(\pi)$.  By Fatou's lemma,
$$\pi(e^{-\lambda_0t}P^V_t|f_n-f|)\le \liminf_m \pi(e^{-\lambda_0t}P^V_t|f_n-f_m|)= \liminf_m \pi(|f_n-f_m|)=\pi(|f_n-f|).$$
Thus $P^V_tf^\pm(x)\le P^V_t|f|(x)<\infty$ for $\pi$-a.a. $x$, $P^V_tf=P^V_tf^+-P^V_tf^-\in L^1(\pi)$, and $P^V_tf_n\to P^V_tf$ in $L^1(\pi)$.
Hence $e^{-\lambda_0t}P^V_t$, $t\ge0$, is a strongly continuous contraction semigroup on $L^1(\pi)$
satisfying (\ref{eq:posmeas}) for $f$ nonnegative Borel. The invariance (\ref{eq:Pinv}) for $f$ in $L^1(\pi)$ follows. \qed

\begin{lemma}
\label{lemma:Q}
Suppose $\pi$ and $\mu$ are measures with $\mu<<\pi$, and let $\psi=d\mu/d\pi$. Then $e^{-\lambda_0t}P^V_t$, $t\ge0$, is a strongly continuous contraction semigroup on $L^1(\pi)$  iff $P^{V,\psi}_t$, $t\ge0$, is a strongly continuous contraction semigroup on $L^1(\mu)$. Moreover, in this case,(\ref{eq:Pinv}) holds for $f$ in $L^1(\pi)$ iff
\beq{eq:Qinv}
\int_XP^{V,\psi}_tf\,d\mu=\int_Xf\,d\mu,\qquad t\ge0,
\eeq
for $f$ in $L^1(\mu)$.  
\end{lemma}

\proof For $u\ge0$ Borel, we have $\norm{u}_{L^1(\mu)}= \norm{u\psi}_{L^1(\pi)}$. Thus $f\in L^1(\mu)$ iff
$f\psi\in L^1(\pi)$. If $f\in L^1(\mu)$, we also have
$$\norm{P^{V,\psi}_tf}_{L^1(\mu)}=\norm{e^{-\lambda_0t}P^V_t(f\psi)}_{L^1(\pi)}$$
and
$$\norm{P^{V,\psi}_tf-f}_{L^1(\mu)}=\norm{e^{-\lambda_0t}P^V_t(f\psi)-f\psi}_{L^1(\pi)}.$$
The result is an immediate consequence of these identities.\qed

A {\em Markov semigroup on $L^1(\mu)$} is a strongly continuous contraction semigroup $Q_t$, $t\ge0$,
on $L^1(\mu)$ satisfying $Q_tf\ge0$ a.s. $\mu$ for $f\ge0$ a.s. $\mu$ and $Q_t1=1$ a.s. $\mu$. 

\begin{lemma}
\label{lemma:QMarkov}
Suppose $\pi$ and $\mu$ are measures with $\mu<<\pi$, and let $\psi=d\mu/d\pi$. Suppose $\pi$ is a ground measure and  $\psi$ is a ground state relative to $\mu$. Then  $P^{V,\psi}_t$, $t\ge0$, is a Markov semigroup on $L^1(\mu)$.
\end{lemma}

\proof 
 $P^{V,\psi}_tf\ge0$  a.s. $\mu$ whenever $f\ge0$ a.s. $\mu$ is clear.  If $\psi$ is a ground state relative to $\mu$, then $P^{V,\psi}_t1=1$ a.s. $\mu$. Thus $P^{V,\psi}_t$, $t\ge0$, is a Markov semigroup on $L^1(\mu)$. \qed

\section{Entropy}
\label{sec:ent}

For $\mu,\pi$ in $M(X)$, the {\em entropy}  of $\mu$ relative to $\pi$ is
$$H(\mu,\pi)\equiv\sup_f\left(\int_Xf\,d\mu-\log\int_Xe^f\,d\pi\right)$$
where the supremum is over $f$ in $C(X)$. 

\begin{lemma} $H(\mu,\pi)\ge0$ is finite iff $\mu<<\pi$ and $\psi\log \psi$  is in $L^1(\pi)$, where $\psi=d\mu/d\pi$,
in which case
$$H(\mu,\pi)=\int_X \psi\log \psi\,d\pi.$$
Moreover $H$ is lower-semicontinuous and convex separately in each of $\mu$ and $\pi$.
\end{lemma}

\proof
The lower-semicontinuity and convexity follow from the definition of $H$ as a supremum of continuous convex functions, in each variable 
$\pi$, $\mu$ separately.

Suppose $H\equiv H(\mu,\pi)<\infty$; then
\beq{eq:H}
\int_Xf\,d\mu-\log\int_Xe^f\,d\pi\le H
\eeq
for $f$ in $C(X)$. The class of Borel functions $f$ for which (\ref{eq:H}) holds is closed under bounded convergence.
Thus (\ref{eq:H}) holds for all $f\in B(X)$. Inserting $f=r1_A$, where $\pi(A)=0$, we obtain
$$r\mu(A)\le r\mu(A)-\log(\pi(A^c))\le H.$$
Let $r\to\infty$ to conclude $\mu<<\pi$. Since $\psi=d\mu/d\pi\in L^1(\pi)$, let $0\le f_n\in C(X)$ with $f_n\to \psi$
in $L^1(\pi)$. By passing to a subsequence, assume $f_n\to \psi$ a.s. $\pi$.
Insert $f=\log(f_n+\epsilon)$ into (\ref{eq:H}) to yield
$$\int_X\log(f_n+\epsilon)\,d\mu-\log\int_X(f_n+\epsilon)\,d\pi\le H.$$
Let $n\to\infty$ followed by $\epsilon\to0$. Since $f_n\to\psi$ in $L^1(\pi)$, 
$$\log\int_X(f_n+\epsilon)\,d\pi\to\log\int_X\psi\,d\pi.$$
Thus, by Fatou's lemma (twice),
$$\int_X \psi\log\psi\,d\pi-\log\int_X\psi\,d\pi\le H.$$
Since $\psi\log^-\psi$ is bounded, this establishes $\psi\log\psi$ in $L^1(\pi)$ and $\int_X\psi\log\psi\,d\pi\le H$.

Conversely, suppose $\psi=d\mu/d\pi$ exists and $\psi\log \psi\in L^1(\pi)$. By Jensen's inequality,
$$\int_X f\,d\mu\le \log\int_X e^f\,d\mu$$
for $f$ bounded Borel. Replace $f$ by $f-\log(\psi\wedge n+\epsilon)$ to get
$$\int_Xf\,d\mu-\log\int_X\left(\frac{e^f\psi}{\psi\wedge n+\epsilon}\right)\,d\pi\le \int_X \psi\log (\psi\wedge n+\epsilon)\,d\pi.$$
Let $\epsilon\to0$ followed by $n\to\infty$ obtaining
$$\int_Xf\,d\mu-\log\int_X e^f\,d\pi\le \int_X \psi\log \psi\,d\pi.$$
Now maximize over $f$ in $C(X)$ to conclude $H(\mu,\pi)\le \int_X \psi\log \psi\,d\pi$.
\qed

\section{Proofs of the theorems}
\label{sec:thm}

{\em Proof of Theorem \ref{theorem:triple}.}
Assume $\pi$ is a ground measure and $\psi$ is a ground state relative to $\mu$. By Lemma \ref{lemma:QMarkov}, $P^{V,\psi}_t$, $t\ge0$, is a Markov semigroup on $L^1(\mu)$. Suppose  $\log u\in L^1(\mu)$. Then 
$P^{V,\psi}_t(\log u)$ is in $L^1(\mu)$, hence there is a set $N$ with $\mu(N)=0$ and $P^{V,\psi}_t(|\log u|)(x)<\infty$ and 
$P^{V,\psi}_t1(x)=1$ for $x\not\in N$. Jensen's inequality applied to the integral 
$f\mapsto (P^{V,\psi}_tf)(x)$ implies 
$$\log\left(\frac{P^{V,\psi}_tu}{u}\right)(x)\ge P^{V,\psi}_t(\log u)(x)-(\log u)(x),\qquad x\not\in N.$$
Thus the negative part of
\beq{eq:integrable}
\log\left(\frac{P^{V,\psi}_tu}{u}\right)
\eeq
is in $L^1(\mu)$ and 
$$\int_X\log\left(\frac{P^{V,\psi}_tu}{u}\right)\,d\mu\ge \int_X \left(P^{V,\psi}_t(\log u)-\log u\right)\,d\mu=0.$$
By Lemma \ref{lemma:nonneg}, this implies $\mu$ is an equilibrium measure, establishing the first claim.

Assume $\pi$ is a ground measure and $\mu$ is an equilibrium measure. Then $e^{-\lambda_0t}P^V_t$, $t\ge0$, is a strongly continuous contraction semigroup on $L^1(\pi)$, hence $P^{V,\psi}_t$, $t\ge0$, is a strongly continuous contraction semigroup on $L^1(\mu)$. Since $1\in L^1(\mu)$, $P^{V,\psi}_t1\in L^1(\mu)$ hence $\log^+(P^{V,\psi}_t1)\in L^1(\mu)$. By Lemma \ref{lemma:nonneg},  $\log^-(P^{V,\psi}_t1)\in L^1(\mu)$ hence
$\log(P^{V,\psi}_t1)\in L^1(\mu)$.
By Jensen's inequality, (\ref{eq:Qinv}), and Lemma \ref{lemma:nonneg},
$$0=\log(\mu(1))=\log\left(\int_X P^{V,\psi}_t1\,d\mu\right)\ge \int_X \log(P^{V,\psi}_t1)\,d\mu\ge0.$$
Since $\log$ is strictly concave, this can only happen if $P^{V,\psi}_t1$ is $\mu$ a.s. constant. By (\ref{eq:Qinv}), the constant
is $1$. Since $\psi>0$ a.s. $\mu$ is immediate, this establishes the second claim.

Assume $\mu$ is an equilibrium measure and $\psi$ is a ground state relative to $\mu$.
Then $P^{V,\psi}_t1=1$ a.s. $\mu$, hence for $u\in B^+(X)$,
$$\frac{\inf u}{\sup u}\le \frac{P^{V,\psi}_tu}{u}\le\frac{\sup u}{\inf u},\qquad a.s. \mu,$$
hence (\ref{eq:integrable})  is in $L^\infty(\mu)$. Thus by (\ref{eq:logQ}), for $f\in B(X)$,
$$\int_X\log\left(\frac{P^{V,\psi}_te^{\epsilon f}}{e^{\epsilon f}}\right)\,d\mu\ge0,\qquad \epsilon>0.$$
Since $P^{V,\psi}_t1=1$ a.s. $\mu$ and  $e^{\epsilon f}=1+\epsilon f+o(\epsilon)$, 
$$\frac1\epsilon\log\left(\frac{P^{V,\psi}_te^{\epsilon f}}{e^{\epsilon f}}\right)\to P^{V,\psi}_tf-f,\qquad \epsilon\to0,$$
uniformly a.s. $\mu$ on $X$. Thus
$$\int_X\left(P^{V,\psi}_tf-f\right)\,d\mu\ge0$$
for $f\in B(X)$. Applying this to $\pm f$ yields
\beq{eq:almost}
\int_Xe^{-\lambda_0t}P^V_t(\psi f)\,d\pi=\int_X\psi f\,d\pi
\eeq
for $f\in B(X)$. Now let $\psi_\delta=(\psi\vee\delta)\wedge(1/\delta)$. Then $\psi/\psi_\delta\to1$ pointwise.
For $f$ in $C(X)$, $f/\psi_\delta\in B(X)$; inserting this into (\ref{eq:almost}) yields
\beq{eq:almost2}
\int_Xe^{-\lambda_0t}P^V_t\left(\frac{\psi f}{\psi_\delta}\right)\,d\pi=\int_X\frac{\psi f}{\psi_\delta}\,d\pi.
\eeq
Since 
$$\frac{|f|\psi}{\psi_\delta}\le  |f|+|f|\psi,\qquad \delta\le1,$$
sending $\delta\to0$ in (\ref{eq:almost2}) yields (\ref{eq:Pinv}). Hence $\pi$ is a ground measure, establishing the third claim.\qed

{\em Proof of Theorem \ref{theorem:pf}.}
Let 
$$M\equiv \sup_{t\ge0} e^{-\lambda_0t}\norm{P^V_t}.$$
By (\ref{eq:log}),
$$\int_X\log\left(\frac{e^{-\lambda_0t}P^{V}_tu}{u}\right)\,d\mu\ge-tI^V(\mu),\qquad u\in C^+(X).$$
Thus for $f\in C(X)$,
$$\int_Xf\,d\mu-\int_X\log\left(e^{-\lambda_0t}P^V_te^f\right)\,d\mu\le tI^V(\mu),\qquad f\in C(X).$$
By Jensen's inequality,
$$\int_Xf\,d\mu-\log\int_X\left(e^{-\lambda_0t}P^V_te^f\right)\,d\mu\le tI^V(\mu),\qquad f\in C(X).$$
Defining
$$Z_t\equiv e^{-\lambda_0t}\mu(P^V_t1)$$
and 
$$\pi_t(f)\equiv\frac{e^{-\lambda_0t}\mu(P^V_tf)}{Z_t}$$
yields
$$\int_Xf\,d\mu-\log\int_Xe^f\,d\pi_t\le tI^V(\mu)+\log Z_t,\qquad f\in C(X).$$
Taking the supremum over all $f$ yields
$$H\left(\mu,\pi_t\right)\le tI^V(\mu)+\log Z_t.$$
Note $Z_t\le M$, $t\ge0$, hence
$$H\left(\mu,\pi_t\right)\le tI^V(\mu)+\log M,\qquad t\ge0.$$
Now set
$$\bar\pi_T(f)\equiv \frac{\int_0^T Z_t\pi_t(f)\,dt}{\int_0^TZ_t\,dt}
= \frac{\int_0^T e^{-\lambda_0t}\mu(P^V_tf)\,dt}{\int_0^TZ_t\,dt},\qquad T>0.$$
Then $\pi_t$ is in $M(X)$ for $t>0$ and $\bar\pi_T$ is in $M(X)$ for $T>0$. 

Now assume $\mu$ is an equilibrium measure; then $I^V(\mu)=0$. By convexity of $H$,
$$H\left(\mu,\bar\pi_T\right)\le \log M,\qquad T>0.$$
By compactness of $M(X)$, select a sequence $T_n\to\infty$ with $\pi_n=\bar\pi_{T_n}$ converging to some $\pi$. By lower-semicontinuity of $H$, we have $H(\mu,\pi)\le \log M$. Thus $\mu<<\pi$ with $\psi=d\mu/d\pi$ satisfying $\psi\log\psi\in L^1(\pi)$. By Lemma \ref{lemma:Iapprox},
$$\log Z_t=\log\mu(e^{-\lambda_0t}P^V_t1)\ge \mu(\log(e^{-\lambda_0t}P^V_t1))\ge0,$$
hence $Z_t\ge1$, $t\ge0$. Since 
$$e^{-\lambda_0t}\mu\left(P^V_t(L+V-\lambda_0)f\right)=\frac{d}{dt} e^{-\lambda_0t}\mu(P^V_tf),\qquad f\in\mathcal D,$$
we have
$$\left|\bar\pi_T((L+V-\lambda_0)f)\right|=\frac{\left|e^{-\lambda_0T}\mu(P^V_Tf)-\mu(f)\right|}{\int_0^TZ_t\,dt}\le\frac{2M\norm{f}}{T}\to0,
\qquad T\to\infty.$$
Thus
$$\pi((L+V-\lambda_0)f)=0,\qquad f\in\mathcal D,$$
which implies
$$e^{-\lambda_0t}\pi(P^V_tf)=\pi(f),\qquad t\ge0,$$
for $f\in\mathcal D$, hence, by density, for $f$ in $C(X)$. Thus $\pi$ is a ground measure.

Clearly $\log\psi\in L^1(\mu)$ and $\psi>0$ a.s. $\mu$. Since $\mu$ is an equilibrium measure and $\pi$ is a ground measure, Theorem \ref{theorem:triple} implies $\psi$ is a ground state relative to $\mu$.\qed


\begin{thebibliography}{1}

\bibitem{A}
S. Aida (1998)
``Uniform Positivity Improving Property, Sobolev Inequalities, and Spectral Gaps''
{\em J. Functional Analysis} {\bf 158}, 152--185.

\bibitem{B}
D. Bakry (2004)
``Functional inequalities for Markov semigroups,''
{\em Probability measures on groups,} 
Probability measures on groups, Mumbai, India. 

\bibitem{DS}
J.-D. Deuschel and D. W. Stroock (1984)
``Large Deviations,''
{\em Pure and Applied Mathematics Series} {\bf 137}, Academic Press.

\bibitem{DV}
M. D. Donsker and S. R. S. Varadhan (1975)
``On a variational formula for the principal eigenvalue for operators with maximum principle,''
{\em Proceedings of the National Academy of Sciences USA} {\bf 72}, 780--783
(\url{http://www.pnas.org/content/72/3/780.full.pdf}).

\bibitem{DV1}
M. D. Donsker and S. R. S. Varadhan (1975)
``Asymptotic evaluation of certain Markov process expectations for large time, I,'''
{\em Communications on Pure and Applied Mathematics} {\bf XXVIII}, 1--47.


\bibitem{E}
R. S. Ellis (1985)
``Entropy, Large Deviations, and Statistical Mechanics,''
{\em Grundlehren der mathematischen Wissenschaften} {\bf 271}, Springer-Verlag.

\bibitem{F}
S. Friedland (1981)
``Convex spectral functions,''
{\em Linear and Multilinear Algebra} {\bf 9}, 299--316.

\bibitem{FK}
S. Friedland and S. Karlin (1975)
``Some inequalities for the spectral radius of  non-negative matrices and applications,''
{\em Duke Mathematical Journal} {\bf 42}, 459--490.

\bibitem{Fr}
G. Frobenius  (1908)
``\"Uber Matrizen aus positiven Elementen,''
{\em S.-B. Preuss. Akad. wiss. Berlin}, 471-476.


\bibitem{GW}
F. Gong and L. M. Wu (2006)
``Spectral Gap of Positive Operators and Applications,''
{\em J. Math. Pures Appl.} [{bf 85} 151--191.

\bibitem{G}
L. Gross  (1972)
``Existence and uniqueness of physical ground states,'' 
{\em J. Func. Anal. } {\bf 10}, 52--109.

\bibitem{JM} 
B. Jourdain and F. Malrieu (2008)
``Propagation of chaos and Poincar\'e inequalities for a system of particles interacting through their cdf ,''
{\em  Annals of Applied Probability} {\bf 18}, 1706-1736 

\bibitem{KR} 
M. G. Krein and M. A. Rutman (1948)
``Linear operators leaving invariant a cone in a Banach space,''
{\em Uspehi Mat. Nauk (N.S.)} {\bf 3}, 3-95. Also { \em Amer. Math.Soc. Transl.} {\bf 26} (1950).


\bibitem{LL}
E. H. Lieb and Michael Loss (1997)
``Analysis,''
{\em Graduate Studies in Mathematics} {\bf 14}, American Mathematical Society.

\bibitem{P}
O. Perron (1907)
``Zur Theorie der Matrices,''
{\em Math. Ann. } {\bf 64}, 248-263.

\bibitem{HHS}
H. H. Schaefer (1971)
``Topological Vector Spaces,''
{\em Graduate Texts in Mathematics} {\bf 3}, Springer-Verlag.

\bibitem{HHS1}
H. H. Schaefer (1974)
``Banach Lattices and positive operators,''
{\em Grundlehren der mathematischen Wissenschaften} {\bf 215}, Springer-Verlag.

\bibitem{SI}
B. Simon (1982)
``Schrodinger Semigroups,''
{\em Bulletin AMS} {\bf 7}, 447-526.

\bibitem{S}
S. Sternberg (2011)
``Dynamical Systems,''
\url{http://www.math.harvard.edu/library/sternberg/}.


\bibitem{W}
L. M. Wu (2000)
``Uniformly Integrable Operators and Large Deviations for Markov Processes,''
{\em J. Functional Analysis} {\bf 172}, 301--376.


\end{thebibliography}
\end{document}